\documentclass{amsart}

\usepackage{amsmath,amssymb}
\usepackage{mathrsfs}
\usepackage{mathtools}
\usepackage{stmaryrd}
\usepackage{enumerate}
\usepackage{tikz-cd}
\usepackage[all]{xy}
\usepackage{aliascnt}
\usepackage[colorlinks=true, linkcolor=blue, citecolor=blue]{hyperref}

\usepackage{enumerate}

\usepackage{combelow}

\newaliascnt{maincorollary}{maintheorem}

\aliascntresetthe{maincorollary}

\newtheorem{theorem}{Theorem}

\newaliascnt{lemma}{theorem}

\aliascntresetthe{lemma}

\newaliascnt{corollary}{theorem}

\aliascntresetthe{corollary}

\newaliascnt{proposition}{theorem}

\aliascntresetthe{proposition}

\newaliascnt{conjecture}{theorem}
\newtheorem{conjecture}[conjecture]{Conjecture}
\aliascntresetthe{conjecture}

\newaliascnt{question}{theorem}

\aliascntresetthe{question}

\theoremstyle{definition}

\newaliascnt{definition}{theorem}

\aliascntresetthe{definition}

\newaliascnt{remark}{theorem}

\aliascntresetthe{remark}

\newaliascnt{example}{theorem}

\aliascntresetthe{example}

\newaliascnt{notation}{theorem}

\aliascntresetthe{notation}

\newif\ifhascomments \hascommentstrue
\ifhascomments
  \newcommand{\matt}[1]{{\color{red}[[\ensuremath{\spadesuit\spadesuit\spadesuit} #1]]}}
  \newcommand{\jeremy}[1]{{\color{red}[[\ensuremath{\clubsuit\clubsuit\clubsuit} #1]]}}
\else
  \newcommand{\matt}[1]{}
  \newcommand{\jeremy}[1]{}
\fi

\newcommand{\bZ}{\mathbb{Z}}

\newcommand{\QQ}{\mathbb{Q}}

\newcommand{\cX}{\mathcal{X}}

\newcommand{\bP}{\mathbb{P}}

\DeclareMathOperator{\st}{st}

\DeclareMathOperator{\Vol}{\Vol}

\tikzset{cong/.style={draw=none,edge node={node [sloped, allow upside down, auto=false]{$\cong$}}},
         Isom/.style={above,every to/.append style={edge node={node [sloped, allow upside down, auto=false]{$\sim$}}}}}

\title[{A counter-example to Batyrev's conjecture}]{A counter-example to Batyrev's conjecture on the non-negativity of stringy Hodge numbers}

\author{Matthew Satriano and Jeremy Usatine}

\thanks{MS was partially supported by an NSERC Discovery Grant. JU was partially supported by the Simons Foundation Travel Support for Mathematicians program and NSF DMS-2502347.}

\address{Matthew Satriano, Department of Pure Mathematics, University of Waterloo}
\email{msatriano@uwaterloo.ca}

\address{Jeremy Usatine, Department of Mathematics, Florida State University}
\email{jusatine@fsu.edu}

\begin{document}

\begin{abstract}
Batyrev's conjecture on the non-negativity of stringy Hodge numbers has been a fundamental open problem, guiding and motivating many beautiful mathematical results in motivic integration, mirror symmetry, and the McKay correspondence. Let $M_0$ be the coarse moduli space of rank 2 semistable bundles with trivial determinant over a fixed smooth projective genus 3 curve. Using a formula, obtained by Kiem and Kiem-Li, for the stringy $E$-function of $M_0$, we verify that $M_0 \times\mathbb{P}^1$ is a counter-example to Batyrev's conjecture.
\end{abstract}

\maketitle

\tableofcontents

\numberwithin{theorem}{section}
\numberwithin{lemma}{section}
\numberwithin{corollary}{section}
\numberwithin{proposition}{section}
\numberwithin{conjecture}{section}
\numberwithin{question}{section}
\numberwithin{remark}{section}
\numberwithin{definition}{section}
\numberwithin{example}{section}
\numberwithin{notation}{section}

\section{Introduction}

Motivated by mirror symmetry for singular Calabi--Yau varieties, Batyrev \cite{Batyrev98} introduced the stringy $E$-function $E_{\st}(X;u,v)$; when $E_{\st}(X;u,v)$ is polynomial, he defined its coefficients, with the appropriate signs, to be the stringy Hodge numbers $h^{p,q}_{\st}(X)$. We emphasize that $h^{p,q}_{\st}(X)$ are defined combinatorially through the use of log-resolutions, rather than as the dimensions of cohomology groups. As a result, positivity is not automatic; nonetheless, Batyrev proposed the following conjecture:

\begin{conjecture}[{\cite[Conjecture~3.10]{Batyrev98}}]
\label{conj:Batyrev}
Let $X$ be a projective algebraic variety with at worst Gorenstein canonical singularities. Assume that $E_{\st}(X; u, v)$ is a polynomial. Then all stringy Hodge numbers $h^{p,q}_{\st}(X)$ are non-negative.
\end{conjecture}

This conjecture has remained one of the most fundamental open problem in motivic integration. More broadly, it has motivated the expectation that stringy Hodge numbers should arise from a cohomology theory naturally associated to a singular variety. Over the past several decades, these questions have guided major developments in motivic integration, mirror symmetry, the McKay correspondence, orbifold cohomology, and the use of Deligne--Mumford and Artin stacks in birational geometry.

When $X$ admits a crepant resolution $Y\to X$, \autoref{conj:Batyrev} follows immediately. Indeed, in this case, $h^{p,q}_{\st}(X)=h^{p,q}(Y)$. This special case is part of a more general story relating motivic integration and crepant resolutions. Motivic integration was pioneered by Kontsevich \cite{Kontsevich} in his proof that birational smooth projective Calabi--Yau varieties have equal Hodge numbers. Denef and Loeser \cite{DenefLoeser2002} further developed the theory, showing that the stringy Hodge numbers can be recovered from the volume of the arc space using the so-called Gorenstein measure. Furthermore, they proved a general motivic change of variables formula yielding a motivic form of the McKay correspondence, refining the numerical correspondence conjectured by Reid \cite{Reid} and proved by Batyrev \cite{Batyrev99}.

Yasuda subsequently proved \autoref{conj:Batyrev} for finite quotient singularities by proving a cohomological interpretation of stringy Hodge numbers in this case. When $X$ has finite quotient singularities, Vistoli \cite{Vistoli89} showed the existence of a canonical smooth Deligne--Mumford stack $\pi\colon\cX\longrightarrow X$ with $\pi$ a small resolution (hence crepant). In his beautiful work \cite{Yasuda2004} (further developed in \cite{Yasuda2006,Yasuda2019}), Yasuda introduced motivic integration over Deligne--Mumford stacks and proved a motivic change of variables formula expressing the Gorenstein measure of $X$ as a motivic integral over the space of twisted arcs of $\cX$. Using this formula, he showed that the stringy Hodge numbers of $X$ agree with the orbifold Hodge numbers of $\cX$ in the sense of Chen--Ruan \cite{ChenRuan}. Thus, for finite quotient singularities, the stringy Hodge numbers are again the dimensions of actual cohomology groups, hence non-negative.

General log-terminal varieties, however, need not admit crepant resolutions by smooth Deligne--Mumford stacks (see, e.g., \cite[Example 7.1]{SatrianoUsatine3}). Nonetheless, in \cite[Theorem 1.1]{SatrianoUsatine3}, the authors showed that every log-terminal variety admits a crepant resolution by a smooth Artin stack. Furthermore, in \cite[Theorem 1.3]{SatrianoUsatine6}, we proved that, for a $\QQ$-Gorenstein variety, the property of admitting such a resolution is \emph{equivalent} to having log-terminal singularities. Using such crepant resolutions, we gave the first known description of stringy Hodge numbers in terms of motivically integrating a function that takes only finitely many values, see \cite[Theorem 1.6]{SatrianoUsatine6}. Using this finite sum expression, we gave a cohomological interpretation of stringy Hodge numbers in \cite{HuangSatrianoUsatine7}. Our \emph{stringy cohomology theory} is naturally mixed and so it does not imply \autoref{conj:Batyrev}; in fact, as we show in this paper, one cannot expect a pure cohomology theory since \autoref{conj:Batyrev} is unfortunately false.

Alongside the search for a cohomological interpretation for stringy Hodge numbers, substantial progress on \autoref{conj:Batyrev} was made directly. Batyrev \cite{Batyrev98} proved non-negativity for complete Gorenstein toric varieties. Schepers and Veys \cite{SchepersVeys} proved \autoref{conj:Batyrev} for threefolds and for several classes of varieties with isolated singularities. Olano \cite{Olano} proved that the stringy Hodge numbers $h_{\st}^{p,1}$ are nonnegative, obtained further results for threefolds even when $E_{\st}$ is not polynomial, and proved Batyrev's conjecture for a class of fourfolds. Other papers on Batyrev's conjecture and related topics include \cite{ChenZuo2025, MustataPayne, SchepersVeys2009, Schepers2009HardLefschetz, Schepers2006ADE}.

There were nevertheless indications that non-negativity was subtle. Schepers and Veys exhibited a projective threefold for which the power-series expansion of the nonpolynomial rational function $E_{\st}(X;u,v)$ has a formally negative value of $h_{\st}^{3,3}$. This was not a counter-example to \autoref{conj:Batyrev} since in their example $E_{\st}(X;u,v)$ was not a polynomial. It did, however, demonstrate that non-negativity is not a formal consequence of the discrepancy formula.

We present the following counter-example to \autoref{conj:Batyrev}, which was discovered with the assistance of ChatGPT.

\begin{theorem}\label{thm:counter-ex}
Let $C$ be a smooth complex projective curve of genus $3$ and let $M_0$ be the coarse moduli space of rank $2$ semistable bundles over $C$ with trivial determinant. Then 
\[
X=M_0\times \bP^1
\]
is a $7$-dimensional projective variety with Gorenstein terminal singularities such that $E_{\st}(X;u,v)$ is a polynomial, but the stringy Hodge number
\[
h_{\st}^{2,5}(X)=-1
\]
is negative. In particular, this is a counter-example to \autoref{conj:Batyrev}.
\end{theorem}

\section{Proof of \autoref{thm:counter-ex}}
The variety $M_0$ has terminal singularities, as shown in \cite[Corollary 5.4]{KiemLi}. Furthermore, it is Gorenstein by the first paragraph of \cite[\S2]{Kiem03}. Since $\dim M_0=3g-3=6$, we see $X$ is $7$-dimensional with Gorenstein terminal singularities. By \cite[Theorem 6.1]{KiemLi}, we see
\[
\begin{aligned}
E_{\st}(M_0;u,v)
={}&
\frac{
 (1-u^2v)^3(1-uv^2)^3
 -(uv)^4(1-u)^3(1-v)^3
}{
 (1-uv)(1-(uv)^2)
}
\\
&-\frac{(uv)^2}{2}
\left(
 \frac{(1-u)^3(1-v)^3}{1-uv}
 -
 \frac{(1+u)^3(1+v)^3}{1+uv}
\right).
\end{aligned}
\]
Thus,
\[
E_{\st}(M_0;u,v)=\frac{P(u,v)}{1+uv},
\]
where
\[
\begin{aligned}
P(u,v) ={} &
u^7 v^7 + u^6 (2 v - 3) v^5 - u^5 (3 v^4 - 3 v^3 + 3 v^2 - 3 v - 1) v^2 - 3 u^4 (v^3 - 4 v^2 - v - 1) v^2 \\
&+ 3 u^3 (v^3 + v^2 + 4 v - 1) v^2 + u^2 (v^4 + 3 v^3 - 3 v^2 + 3 v - 3) v - u (3 v - 2) v + 1.
\end{aligned}
\]
As a result,
\[
E_{\st}(X;u,v)=E_{\st}(M_0;u,v)E(\bP^1;u,v)=P(u,v)\in\mathbb \bZ[u,v],
\]
i.e., $E_{\st}(X;u,v)$ is polynomial.

Lastly, to compute $h^{2,5}_{\st}(X)$, note that the $u^2v^5$-coefficient of $E_{\st}(X;u,v)$ is $1$, hence, by Batyrev's definition,
\[
h_{\mathrm{st}}^{2,5}(X)=(-1)^{2+5}\cdot 1=-1.
\]
This completes the proof.

\bibliographystyle{alpha}
\bibliography{Batyrev-counterex}

\end{document}